# Analytic renormalization of multiple zeta functions. Geometry and combinatorics of the generalized Euler reflection formula for MZV

Andrei Vieru

andreivieru007@gmail.com

**Abstract**

The renormalization of MZV was until now carried out by algebraic means. We show that renormalization in general, of the multiple zeta functions in particular, is more than mere convention. We show that simple calculus methods allow us to compute the renormalized values of multiple zeta functions in any dimension for arguments of the form (1,…,1), where the series do not converge. These values happen to be the coefficients of the asymptotic expansion of the inverse Gamma function. We focus on the geometric interpretation of these values, and on the combinatorics their closed form encodes, which happen to match the combinatorics of the generalized Euler reflection formula discovered by Michael E. Hoffman, which in turn is a kind of analogue of the Cayley-Hamilton theorem for matrices. By means of one single limit formula, we define a function on the positive open half-line which takes exactly the values of the Riemann zeta function, with the additional advantage that it equals the Euler constant when the argument is 1.

**1. Introduction**

Multiple zeta functions are defined as

$$\zeta(s_1,\ldots,s_k) := \sum_{0<n_k<\ldots<n_1} \frac{1}{n_1^{s_1}\ldots n_k^{s_k}} \qquad (1.1)$$

A lot is known about multiple zeta functions and their special values and their various and beautiful properties. The renormalization of these functions, which can be meromorphically extended to complex values and have singularities in many points[1], is closely related to complicate algebraic structures. Hopf algebras and quasi-shuffle relations are cited as soon as renormalization is spoken about. The need of preserving algebraic structures has given birth to the beautiful term "algebraic continuation".

Our start point is not rooted in algebra but in very simple Eulerian calculus.

So instead of "algebraic continuation" we shall speak about a completely different technique that we'll call *extension by discontinuity*, both by analogy and contrast with the well-known "extension by continuity".

---

[1] see [1]

We first redefine the multiple zeta functions, starting with the Riemann zeta function itself, completely giving up series and substituting to them limit definitions.

As we already showed in another paper, Euler might have defined the celebrated function on positive real numbers ≥ 1 not as the series

$$\zeta(s) = \sum_{n=1}^{\infty} \frac{1}{n^s} \tag{1.2}$$

but as the limit

$$\zeta(s) = \lim_{k \to \infty} \left\{ \sum_{n=1}^{k} \frac{1}{n^s} - \sum_{m=k+1}^{k^2} \frac{1}{m^s} \right\} \tag{1.3}$$

The advantage of this definition consists in the fact that **(1.3)**, unlike **(1.2)**, converges for *s* = 1, taking exactly the same values as in **(1.2)** for *s* >1 (when **(1.2)** converges, the subtracted sum in **(1.3)** tends to 0 as *k*→∞).

According to **(1.3)**, $\zeta(1) = \gamma$ (where γ is the Euler constant), which we consider to be the normal value, or if one prefers, the "renormalized value" of the Riemann zeta function for *s* = 1.

The geometric interpretation of the case *s* = 1 is easy to understand if one presents the logarithm of *k* in one of its integral forms, namely $\ln k = \int_{k}^{k^2} \frac{1}{x} dx$ , and then replaces it by the sum at integers values.

This analytic definitions makes it easier to understand why Euler's constant behaves so unambiguously as ζ(1), in particular, in formulas where the values of the Riemann zeta function at integer arguments *k* show up along with Harmonic numbers of order *k*. For example

$$\ln \Gamma\left(-n + (-1)^n x\right) = -\ln x - \sum_{j=1}^{n} \ln j + \sum_{k=1}^{\infty} (-1)^{(n+1)k} \frac{\zeta(k) + (-1)^k H_n^{(k)}}{k} x^k \tag{1.4}$$

(for *n* ≥ 1 and when 0 < *x* ≤ 1/2)

We also have the following formula, which involves the Riemann zeta function at odd arguments[2]:

$$\ln n = \sum_{k=1}^{\infty} \left( \frac{1}{2kn^{2k}} \frac{1}{2^{2k}} - \frac{2\zeta(2k-1) - H_{n-1}^{(2k-1)} - H_n^{(2k-1)}}{2k-1} \frac{1}{2^{2k-1}} \right) \tag{1.5}$$

(which, using Harmonic numbers of fractional argument, holds as well for fractional *n* ≥ 1)

The Euler constant behaves as ζ(1), as well, in recursion rules, as for instance in **(3.6)** (see bellow)

---

[2] In both **(1.4)** and **(1.5)** ζ(1) is taken to be the Euler constant.

## 2. An elementary example of renormalization through "extension by discontinuity" concerning the multiple zeta function in two variables

We shall place ourself in a purely Eulerian perspective: Euler defined the zeta function in one and two *real* variables We shall first consider, as a simple example, the case $k = 2$ with $s_1 = s_2 = 1$, where there is a singularity.

For the zeta function in two real variables 
$$\zeta(s_1, s_2) := \sum_{0 < n_2 < n_1} \frac{1}{n_1^{s_1} n_2^{s_2}} \quad (2.1)$$

we propose an alternative limit definition for real arguments $s_1 \geq 1$, $s_2 \geq 1$ according to which

$\zeta(s_1, s_2)$ is, when $m \to \infty$, the limit of

$$\sum_{n_2=1}^{m} \sum_{n_1=1}^{m} \frac{1}{n_1^{s_1} n_2^{s_2}} - \sum_{n_2=1}^{m} \sum_{n_1=m+1}^{m^2} \frac{1}{n_1^{s_1} n_2^{s_2}} - \left( \sum_{n_2=m+1}^{m^2} \sum_{n_1=1}^{m} \frac{1}{n_1^{s_1} n_2^{s_2}} - \sum_{n_2=m+1}^{m^2} \sum_{n_1=m+1}^{m^2} \frac{1}{n_1^{s_1} n_2^{s_2}} \right)$$

$$- H_{m^2}^{(2)} - \sum_{0 < n_1 < n_2 \leq m} \frac{1}{n_1^{s_1} n_2^{s_2}} \quad (2.2)$$

For the points $(s_1, s_2)$ in which **(2.1)** converges, **(2.2)** is equivalent to **(2.1)**: when $(s_1, s_2)$ converges only the first double sum does not vanish. When the two last summands are subtracted from the first double sum and when $m \to \infty$, we get exactly the values of $\zeta(s_1, s_2)$. The arguments of the summands in the first double sum lie in a square with side $m$, respectively $m^2$ if we take into account the three other double sums whose integer "arguments[3]" run within two rectangles with sides $m$ and $m^2-(m+1)$ and within a square with side $m^2-(m+1)$. The subtracted last sum in **(2.2)** corresponds to the half part of the square that lies above the main bisectrix. Its subtraction rules out the summands with $n_1 > n_2$, while the subtraction of the Harmonic number of order 2 and argument $m^2$ eliminates the case $n_1 = n_2$ which is also ruled out by the definition **(2.1)**.

In the special case when $s_1 = s_2 = 1$, for which **(2.1)** does not converge, it is quite easy to show that **(2.2)** reduces to[4]

---

[3] of course, not to be confused with variables of the MZ functions, which are exponents. Here we are speaking about the $n_i$

[4] for this case, in **(2.2)**, using **(1.3)** and passing to the limit in one direction (i.e. with the respect to the first variable), the difference between the two first double sums may be viewed as $\gamma(1+1/2+\ldots+1/m)$, while the difference in the parentheses may be viewed as $\gamma \ln m$. Passing to the limit in the other direction – i.e. with respect to the second variable – one gets $\gamma^2$ for all four double sums in **(2.2)**. We believe that the arithmetic and geometric sense of **(2.2)** is clear enough…

$$\zeta(1,1) := \frac{1}{2} \lim_{m \to \infty} \left[ \sum_{n_2=1}^{m} \sum_{n_1=1}^{m} \frac{1}{n_1 n_2} - 2 \sum_{n_2=m+1}^{m^2} \sum_{n_1=1}^{m} \frac{1}{n_1 n_2} + \sum_{n_2=m+1}^{m^2} \sum_{n_1=m+1}^{m^2} \frac{1}{n_1 n_2} - H_{m^2}^{(2)} \right]$$

$$= \frac{\gamma^2 - \zeta(2)}{2} \tag{2.3}$$

Since $\zeta(1) = \gamma$, the constant in the RHS of **(2.3)** satisfies the "quasi-shuffle relation"

$$\zeta(s_1)\zeta(s_2) = \zeta(s_1, s_2) + \zeta(s_2, s_1) + \zeta(s_1 + s_2) \tag{2.4}$$

which is also known as the "Euler's reflection formula":

## 3. Renormalized values of Multiple zeta functions of higher dimension

One can see that the method of "extension by discontinuity" can be easily carried out for points of the form (1,1,…,1) with arbitrarily large number of unit arguments.
The standard definition of the multiple zeta function in three arguments reads:

$$\zeta(s_1, s_2, s_3) := \sum_{0 < n_3 < n_2 < n_1} \frac{1}{n_1^{s_1} n_2^{s_2} n_3^{s_3}} \tag{3.1}$$

For three arguments we get as an alternative definition the multiple zeta function a quite long limit formula (where $m \to \infty$) which should obviously be divided by 6:

$$\sum_{n_3=1}^{m}\sum_{n_2=1}^{m}\sum_{n_1=1}^{m} \frac{1}{n_1^{s_1} n_2^{s_2} n_3^{s_3}} - \sum_{n_3=1}^{m}\sum_{n_2=1}^{m}\sum_{n_1=m+1}^{m^2} \frac{1}{n_1^{s_1} n_2^{s_2} n_3^{s_3}} - \left( \sum_{n_3=1}^{m}\sum_{n_2=m+1}^{m^2}\sum_{n_1=1}^{m} \frac{1}{n_1^{s_1} n_2^{s_2} n_3^{s_3}} - \sum_{n_3=1}^{m}\sum_{n_2=m+1}^{m^2}\sum_{n_1=m+1}^{m^2} \frac{1}{n_1^{s_1} n_2^{s_2} n_3^{s_3}} \right)$$

$$- \left[ \sum_{n_3=m+1}^{m^2}\sum_{n_2=1}^{m}\sum_{n_1=1}^{m} \frac{1}{n_1^{s_1} n_2^{s_2} n_3^{s_3}} - \sum_{n_3=m+1}^{m^2}\sum_{n_2=1}^{m}\sum_{n_1=m+1}^{m^2} \frac{1}{n_1^{s_1} n_2^{s_2} n_3^{s_3}} - \left( \sum_{n_3=m+1}^{m^2}\sum_{n_2=m+1}^{m^2}\sum_{n_1=1}^{m} \frac{1}{n_1^{s_1} n_2^{s_2} n_3^{s_3}} - \sum_{n_3=1}^{m^2}\sum_{n_2=m+1}^{m^2}\sum_{n_1=m+1}^{m^2} \frac{1}{n_1^{s_1} n_2^{s_2} n_3^{s_3}} \right) \right]$$

$$-3 \left( \sum_{n=1}^{m} \frac{1}{n^s} - \sum_{n=m+1}^{m^2} \frac{1}{n^s} \right) H_{m^2}^{(2)} + 2 H_{m^2}^{(3)} \tag{3.2}$$

If one of the the three variables $s_1$, $s_2$, $s_3$ is strictly greater than 1, then all triple sums in **(3.2)** except the first one vanish as $m \to \infty$ and for this reason **(3.2)** is equivalent to **(3.1)**.
The eight first summands are triple sums, which the method of passing to the limit briefly described in note 2 is liable to be applied to. If $s_1 = s_2 = s_3 = 1$, then they yield $\gamma^3$ after one passes to the limit thrice, with respect to each variable, and of course using **(1.3)**. So the arguments under the triple sums run within a cube of side $m$ (with $m \to \infty$), which one subtracts from and/or adds to (if one gets rid of the parentheses) "parallelepipeds" with various combinations of lengths of sides $m$ and $m^2 - (m + 1)$, and a "cube" with side $m^2 - (m + 1)$. So everything is happening within a cube with side $m^2$ (with $m \to \infty$).

(We write write "parallelepipeds" and "cube" with quotation marks because we do not integrate but just sum.)

We have to subtract thrice
$$\left(\sum_{n=1}^{m} \frac{1}{n^s} - \sum_{n=m+1}^{m^2} \frac{1}{n^s}\right) H_{m^2}^{(2)} \qquad (3.3)$$

because we have to rule out, according to the standard definition **(3.1)**, the cases $n_1 = n_2$, $n_2 = n_3$ and $n_1 = n_3$; according to the same standard definition, we must also rule out the case $n_1 = n_2 = n_3$ but *only once*. And it was already subtracted thrice while ruling out the cases $n_1 = n_2$, $n_2 = n_3$ and $n_1 = n_3$, since each of them in fact include the latter.

So we have to add $H_{m^2}^{(3)}$ twice. Now we have to divide the whole stuff by 6, for we have, according to the standard definition, one single admitted order $n_3 < n_2 < n_1$ out of six possible permutations of this order. After a simple computation, knowing that the limit of **(3.3)** when $m \to \infty$ is $\gamma\zeta(2)$, we get as a limit of **(3.2)** when $m \to \infty$ the following value:

$$\zeta(1,1,1) = \frac{\gamma^3 - 3\gamma\zeta(2) + 2\zeta(3)}{6} \qquad (3.4)$$

For dimension 4, the alternative limit definition of the multiple zeta function is to heavy and we give up writing it down. It is not difficult to guess how to write the quadruple sums in order to get all of them (save the first) vanishing when at least one argument $s_i$ ($1 \leq i \leq 4$) is greater than 1. Care has to be shown when we subtract and/or re-add summands which rule out cases of equality between the running integers under summation. And once again the whole stuff has to be eventually divided by 24, for the standard definition of the multiple zeta function prescribes a fixed order, which is only one of the 24 possible permutations. (Once again we are allowed to take into account symmetry, since the arguments all equal 1.) Finally we get the value

$$\zeta(1,1,1,1) = \frac{\gamma^4 - 6\gamma^2\zeta(2) + 3\zeta(2)^2 + 8\gamma\zeta(3) - 6\zeta(4)}{24} \qquad (3.5)$$

By the way, we might restrain from using further the $\gamma$ symbol and replace it everywhere by $\zeta(1)$. We've got the fifth coefficient in the well-known expansion of $1/\Gamma(z)$. (Let us consider 1 as being the renormalized value of the multiple zeta function in 0 variables, and let us assign it the rank 0.) It can be shown that for higher dimensions we shall still get "renormalized" values the coefficients of $1/\Gamma(z)$

It is well-known that the coefficients of $1/\Gamma(z)$ satisfy the following recursion rule:

$$a_0 = 1$$
$$a_1 = \gamma$$
$$2a_2 = \gamma a_1 - \zeta(2)$$

…

$$na_n = \sum_{k=1}^{n}(-1)^{k+1}a_{n-k}\zeta(k) \qquad \text{(valid for } n \geq 1, \text{ and taking } \varsigma(1) = \gamma) \tag{3.6}$$

The "renormalized" values $\varsigma(1,\ldots,1)$, which perfectly match the coefficients of $1/\Gamma(x)$, may be considered, formally, as polynomials in "pseudo-variables" of the form $\zeta(k)$. Replacing in these forms $\zeta(k)$ by $tr(A^k)$, and in particular $\gamma = \zeta(1)$ by $tr(A)$, one gets the RHS of the well-known infinite collection of identities for 1×1, 2×2, 3×3, 4×4 (and so on) matrix determinants in terms of traces. Thus, $\det(A) = tr(A)$ in dimension 1×1. For the next dimensions, we have:

$$\det(A) = \frac{(tr(A))^2 - tr(A^2)}{2} \qquad \text{which has to be compared with } \textbf{(2.3)}$$

$$\det(A) = \frac{(tr(A))^3 - 3tr(A)tr(A^2) + 2tr(A^3)}{6} \qquad \text{which has to be compared with } \textbf{(3.4)}$$

$$\det(A) = \frac{(tr(A))^4 - 6tr(A^2)(tr(A))^2 + 3(tr(A^2))^2 + 8tr(A^3)tr(A) - 6tr(A^4)}{24} \qquad \text{(to be compared with } \textbf{(3.5)})$$

and so on…

Actually, we don't know whether the integer coefficients in the coefficients of the Taylor series of $1/\Gamma(z)$ were studied or not as such. But they were indeed studied in the context of the Newton identities and Cayley-Hamilton theorem. If taken unsigned, as they appear[5] either in the Taylor series of $\Gamma(1 + z)$ or in the Laurent series of $\Gamma(z)$ near 0, they are nothing else but the "array of multinomial numbers", which can be found in OEIS [A102189](#) or [A036039](#). At the moment we are writing this paper, nothing is said there about the Gamma function and the related series or — and that's indeed more important — about the generating recursion rule. One can see a lot of links in those webpages of OEIS that my provide possible connections with numerous domains of research or knowledge. The most ready to hand connexion with elementary Algebra is suggested by the fact that they appear in cycle indexes of permutation groups (symmetric groups for example). They are closely related to partitions and to the partition function[6] ([A000041](#) in [OEIS](#)), which are known to occur, for example, in Group representation theory.

---

[5] the coefficients themselves appear signed, but we are speaking about the integer coefficients in the numerator of these coefficients

[6] the number of summands in the numerators of this values is given by the partition function.

## 4. A geometric-combinatorial interpretation of the signed integer coefficients in the numerators of the renormalized values of ζ(1,…, 1)

As one easily may imagine this interpretation concerns the diagonal line, plans, and hyperplans of the hypercube that contain the points $(1, 1,…, 1)$ and $(m, m,…, m)$. (These hyperplans, plans and line will be called *main diagonals*). In dimension 3, we have integer coefficients 1, –3 and 2 (see the LHS of **(3.4.)**)

In the figure below, one can see the upper facet of a cube, which is its own 3-dimensional main diagonal, 3 plans corresponding, respectively to the equations $x = y$, $x = z$ and $y = z$, and a line $x = y = z$ where the three plans meet together[7].

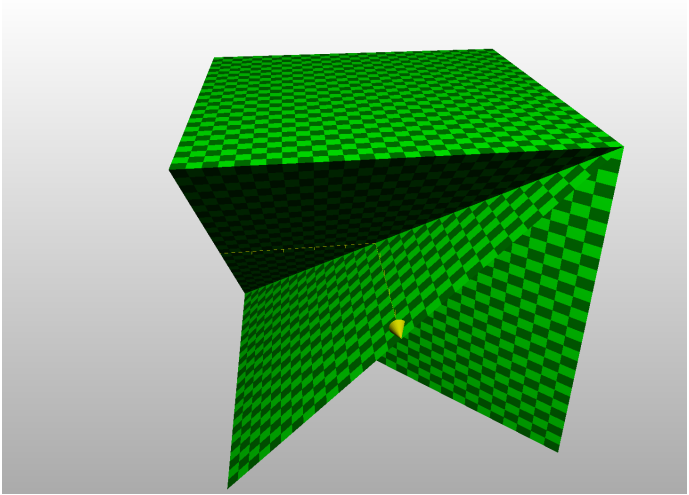

The integers 1, –3 and 2 correspond to the subtraction *only once* from the cube of its lower-dimensional diagonals that contain $(1, 1, 1)$ and $(m, m, m)$. As we pointed out in section **2.**, after we subtracted the three plans, and therefore subtracted thrice the line $x = y = z$, we have to readd it twice. The combinatorics of the integer coefficients in the numerator of the Euler-McLaurin series of $1/\Gamma(z)$ are exactly the combinatorics of the subtraction-readdition of lower-dimensional main diagonals from a finite or infinite hypercube, *when one wants the final result to be a subtraction of each main diagonal only once.*

Clearly, all diagonals are determined by a combination of equalities. Thus, in an *n*-dimensional cube the main diagonal line is determined by the equality between all variables, while the cube itself is determined by the equality of each variable to itself. These combinatorics reduce to the combinatorics of the intersection, union and subtraction (∪, ∩, \ or, accordingly, if considered in the frame of propositional calculus ∧,∨,\) of these equalities which in turn imply other combinations of equalities between the variables as special cases. It is easy to see that these combinatorics are closely related to partitions and to decorated rooted trees without sidebranchings.

---

[7] here the plans are of course restricted to the magnitude of cube

## 5. Quasi-shuffle relations over multiple zeta functions of all dimensions

Even if one knows nothing about partitions, partition numbers, or decorated rooted trees without sidebranchings the integer coefficients in the value expressed by **(3.4)** suggests the following relation, which is indeed satisfied by the special value of $\zeta(1,1,1)$ itself, and is in fact a well-known[8] generalization of the Euler's reflection formula **(2.4)**

$$\zeta(s_1)\zeta(s_2)\zeta(s_3) = \zeta(s_1,s_2,s_3) + \zeta(s_1,s_3,s_2) + \zeta(s_2,s_3,s_1) + \zeta(s_2,s_1,s_3) + \zeta(s_3,s_1,s_2) + \zeta(s_3,s_2,s_1)$$
$$+\zeta(s_1)\zeta(s_2+s_3) + \zeta(s_2)\zeta(s_1+s_3) + \zeta(s_3)\zeta(s_1+s_2) - 2\zeta(s_1+s_2+s_3) \qquad (5.1)$$

Obviously, the unsigned terms of **(5.1)** correspond exactly to a complete collection of partitions or, simpler to say in this context, of to a complete collection of decorated rooted trees without sidebranchings, as do the summands in the numerator of **(3.4)**.

The renormalized value at $\zeta(1,1,1,1)$ encodes a "quasi-shuffle" relation of the same type. Any summand in the numerator of this value – see **(3.5.)** – suggests the content of the suitable corresponding term (or group of terms) in the RHS of the "quasi-shuffle" relation: $\zeta(k)$ always corresponds to an equality (ruled out or restored, it depends on the sign) between exactly $k$ terms in our "hyper-cubic model" of the alternative limit definition of the Multiple Zeta functions. This equalities (or better to say complete collection of systems of equalities, for we may have, say, $s_1=s_4$ and $s_2=s_3$ or any other conceivable combination of equalities) translate in a complete collection of decorated rooted trees without sidebranchings in the RHS of the "quasi-shuffle" relation.

In the general *normal* case, these relations are completely described in Theorem 2.2 in [3]. They still hold for renormalized values at singularity points of the form $(1,…,1)$. As far as we know, nothing was never said about the elementary fact that these combinatorics perfectly match the combinatorics of the coefficients of $1/\Gamma(z)$, and us such are completely described by the recursion rule **(3.6)**

## 6. "Extension by discontinuity": analytic renormalization of some simple integrals

One can redefine the Gamma function as the following limit:

$$\Gamma(x) = \lim_{y \to 0}\left[\int_y^\infty t^{x-1}e^{-t}dt - \int_{y^2}^y t^{x-1}e^{-t}dt\right]$$

---

[8] see [2] and [3]

This definition is equivalent to Euler's definition in its well-known integral form with only one significant difference: Euler's integral diverges when $x = 0$, while our limit converges. According to our limit definition $\Gamma(0) = -\gamma$, which is the renormalized, or – as we say – the *normal* value of $\Gamma(0)$.

Define for $x \geq 1$ the function $\quad \varpi(x) = \int_0^{\pi/2} (\tan u)^{\frac{1}{x}} \, du$

Oh sorry, when $x = 1$ the integral diverges… So redefine the same function as a limit, where $y$ is supposed to be nonnegative:

$$\varpi(x) = \lim_{y \to 0} \left[ \int_0^{\pi/2-y} (\tan u)^{\frac{1}{x}} \, du - \int_{\pi/2-y}^{\pi/2-y^2} (\tan u)^{\frac{1}{x}} \, du \right]$$

then $\quad \varpi(1) = 0$

The cosine integral is defined as $\quad Ci(x) = -\int_x^\infty \frac{\cos t}{t} \, dt \quad$ and it diverges when $x = 0$

If redefined as $\quad Ci(x) = \lim_{y \to 0} \left[ -\int_{x+y}^\infty \frac{\cos t}{t} \, dt + \int_{x+y^2}^{x+y} \frac{\cos t}{t} \, dt \right] \quad$ then $Ci(0) = \gamma$

The last example shows that the rate of growth (or decay) of the integrand may matter as much as the rate of growth (or decay) of the integral itself. This fact is obvious when the variable is written as the lower (or upper) bound of the integral.

For example, define $\quad Ci_2(x) = \int_x^\infty \frac{\cos t}{t^2} \, dt$

The integral diverges when $x = 0$, but one can redefine the function as

$$Ci_2(x) = \lim_{y \to 0} \left[ \int_{x+y}^\infty \frac{\cos t}{t^2} \, dt - \int_{x+\frac{y}{2}}^{x+y} \frac{\cos t}{t^2} \, dt \right] \quad \text{and get } Ci_2(0) = -\pi/2$$

In contrast with the previous examples, in the last one the bounded variable $y$ does not appear anymore quadratically, but the value $-\pi/2$ is, as well as in all other examples, the constant term of the Laurent expansion near the singularity of the function given in its integral form.
We believe it is important to make it clear that choosing the constant terms of Laurent series as renormalized values of functions at their singularities is more than a pure algebraic convention.
Introducing a parameter $r$ (see also [5] for the possible role played by an additional parameter when Riemann zeta function is renormalized at $s=1$), one can write more generally:

$$\Gamma(x) = \lim_{y \to 0} \left[ \int_y^\infty t^{x-1} e^{-t} dt - \int_{ry^2}^y t^{x-1} e^{-t} dt \right] = -\gamma + \ln r$$

$$\varpi(x) = \lim_{y \to 0} \left[ \int_0^{\pi/2-y} (\tan u)^{\frac{1}{x}} du - \int_{\pi/2-y}^{\pi/2-ry^2} (\tan u)^{\frac{1}{x}} du \right] = \ln r$$

$$Ci(x) = \lim_{y \to 0} \left[ -\int_{x+y}^\infty \frac{\cos t}{t} dt + \int_{x+ry^2}^{x+y} \frac{\cos t}{t} dt \right] = \gamma - \ln r$$

However, the natural way of analytic renormalization is to set $r = 1$

## 7. Further "extension by discontinuity" of the Riemann zeta function

The limit definition **(1.3)** can be further extended to divergent $\sum_{k=1}^\infty \frac{1}{k^s}$ series where $0 < s < 1$

Let us consider the integral $\int_1^n \frac{1}{x^s} dx$

Under the presupposed condition $s \neq 1$ the antiderivative of the integrand is $\frac{x^{1-s}}{1-s}$

To get a definition in the style of **(1.3)** one has to solve the equation

$$\frac{(2^y x)^{1-s}}{1-s} = \frac{2x^{1-s}}{1-s} \qquad \text{(where the unknown is } y\text{)} \qquad \text{(7.0)}$$

The solution is obvious: when $y=1/(1-s)$ the equality **(7.0)** becomes an identity, and it will be used as a power of 2 in order to fix the upper bound of the subtracted sum (using if necessary the floor function $\lfloor t \rfloor$)

We state that for strictly positive $s$ smaller than 1, where the Euler-Riemann series diverges,

$$\zeta(s) = \lim_{m \to \infty} \left[ \sum_{n=1}^m \frac{1}{n^s} - \sum_{k=m+1}^{\lfloor 2^{\frac{1}{1-s}} m \rfloor} \frac{1}{k^s} \right]$$

For the elegance, one can define a family of functions $f_s(t)$ via the family of functional equations

$$\int_1^t \frac{dx}{x^s} = \int_t^{f_s(t)} \frac{dx}{x^s} + C \tag{7.1}$$

which, for a given $s > 0$, is an equality supposed to hold for all $t$, and where $C$ is a constant[9], namely the smallest real number in absolute value for which the functional equation **(7.1)** has a monomial solution. For any strictly positive $s$, it is quite easy to solve **(7.1)**. Everyone knows that $f_1(t) = t^2$ is a solution of the functional equation **(7.1)** in which $s$ is set to 1, and $C$ is set to 0. Besides, if $s$ is greater than 0 and smaller than 1, then $f_s(t) = 2^{1/(1-s)}t$ and $C = 1/(1-s)$

Therefore, we can define a function $$\zeta(s) = \lim_{m \to \infty} \left[ \sum_{n=1}^{m} \frac{1}{n^s} - \sum_{k=m+1}^{\lfloor f_s(m) \rfloor} \frac{1}{k^s} \right] \tag{7.2}$$

where $f_s$ is a solution of **(7.1)**. This function takes exactly the values of the Riemann zeta function for any strictly positive $s$ smaller than 1, and equals the Euler constant when $s = 1$.

Let us consider the family of 'limit functional equations'

$$\lim_{t \to \infty} \left\{ \int_1^t \frac{dx}{x^s} - \int_t^{g_s(t)} \frac{dx}{x^s} \right\} = C(s) \tag{7.3}$$

If $s \leq 1$, then $g_s = f_s$

If $s > 1$, then the first integral under the limit converges when $t = \infty$, and one can consider the trivial solutions $g_s(t) = t$. For a given $s$, the constant $C$ is therefore nothing but the evaluation of the first integral between 1 and infinity.

By analogy with the functions $f_s$, when $s > 1$, one can consider as a family of solutions of **(7.3)** $g_s(t) = 2^{1/(1-s)}t$ (or, by the way, any strictly increasing unbounded function).

Replacing $f_s$, by $g_s$ (i.e. using the solutions of the family of functional equations **(7.3)** instead of **(7.1)**) one can view **(7.2)** as a unique definition of the Riemann zeta function valid for the entire strictly positive real domain, including the argument 1.

When $s$ is of the form $n/(n+1)$, there is no need to use the floor function. The smallest value is 1/2 (which belongs to the "critical line"), and we have $f_{0.5}(m) = 2^{1/(1-1/2)}m = 4m$

The convergence is very slow in the strictly decreasing sequence **(7.2)** as $m \to \infty$. For example:

$$\sum_{1}^{15000} \frac{1}{\sqrt{n}} - \sum_{15001}^{60000} \frac{1}{\sqrt{k}} = -1.4542\ldots$$

while $\zeta(1/2) = -1.460355\ldots$ Skipping one summand, we get a strictly increasing sequence with a slightly faster convergence. For example:

$$\sum_{1}^{14999} \frac{1}{\sqrt{n}} - \sum_{15001}^{60000} \frac{1}{\sqrt{k}} = -1.4624\ldots$$

---

[9] i.e. a real number that does not depend on $t$ (but may depend on $s$ if for a given $s$ one wants the equation to have a solution)